\documentclass[
 reprint,
superscriptaddress,
 amsmath,amssymb,prl
 aps
]{revtex4-2}

\usepackage{natbib}
\usepackage{amsmath,amssymb,amsfonts,graphicx,url,algorithm2e}
\usepackage[dvipsnames]{xcolor}
\usepackage{times}
\usepackage{dsfont}
\usepackage{mathtools}
\usepackage[utf8]{inputenc}
\usepackage[T1]{fontenc}
\usepackage{fancyhdr}
\usepackage{subfigure}
\usepackage{enumerate}
\graphicspath{{./},{./figures/}}

\newcommand{\dbar}{d\hspace*{-0.1em}\bar{}\hspace*{0.2em}}

\newcommand{\ud}{d} 

\newcommand{\mR}{\mathbb{R}}

\newcommand{\calP}{\mathcal{P}}

\newcommand{\etaU}{\eta_{\hspace*{-2pt}\;_\mathcal{U}}}
\newcommand{\etaQ}{\eta_{\hspace*{-2pt}\;_\mathcal{Q}}}

\newcommand{\magenta}{\color{black}}
\newcommand{\red}{\color{black}}

\newcommand{\black}{\color{black}}
\newcommand{\purple}{\color{black}}
\usepackage{comment}
\usepackage{graphicx}
\usepackage{dcolumn}
\usepackage{bm}

\begin{document}

\preprint{APS/123-QED}

\title{
Underdamped stochastic thermodynamic engines\\ in contact with a heat bath with arbitrary temperature profile}%

\author{Olga Movilla Miangolarra}
\affiliation{%
Department of Mechanical and Aerospace Engineering, University of California, Irvine, CA 92697, USA
}
\author{Rui Fu}%
\affiliation{%
Department of Mechanical and Aerospace Engineering, University of California, Irvine, CA 92697, USA
}
\author{Amirhossein Taghvaei}
\affiliation{%
Department of Mechanical and Aerospace Engineering, University of California, Irvine, CA 92697, USA
}%
\author{Yongxin Chen}
\affiliation{%
School of Aerospace Engineering, Georgia Institute of Technology, Atlanta, GA 30332, \magenta{USA}
}%
\author{Tryphon T. Georgiou}
\affiliation{%
Department of Mechanical and Aerospace Engineering, University of California, Irvine, CA 92697, USA
}%

\date{\today}

\begin{abstract}
	We study thermodynamic processes in contact with
	a heat bath that may have an arbitrary time-varying periodic temperature profile.
	Within the framework of stochastic thermodynamics, and for models of thermodynamic engines in the idealized case of underdamped particles in the low friction regime {subject to a harmonic potential}, we derive explicit bounds as well as optimal control protocols that draw maximum power and achieve maximum efficiency at any specified level of power.

\end{abstract}

\maketitle

	For more than 200 years, Carnot's idealized concept of a heat engine \cite{carnot1986reflexions} 
	has been the cornerstone of equilibrium thermodynamics {{\cite{adkins1983equilibrium}}}.
	It provided a standard as well as led to the discovery of the absolute temperature scale, entropy,
	and the second law of thermodynamics.
	The sequel, beyond equilibrium and quasi-static operation, proved a challenging one.
	Statistical mechanics and stochastic fluctuation theories have since, and for over a hundred years,
	been advancing our insight into the irreversibility of nonequilibrium thermodynamic transitions.
	In this endeavor, a particular promising model of far-from-equilibrium transitions is that of Stochastic Thermodynamics (also, Stochastic Energetics) {{\cite{sekimoto2010stochastic,seifert2012stochastic,chen2019stochastic,Jarz2011ineq}}} that followed
	at the heels of fluctuation theories \cite{Jarz1996eq,Cohen1995FT,Evans1994FT,Crooks1999FT,Hatano2001FT}.
	
	In the present work we adopt the framework of Stochastic Thermodynamics so as to quantify the amount of power as well as the efficiency that thermodynamic engines can achieve
	when harvesting energy from a heat bath with an arbitrary but periodically fluctuating temperature profile.
	Our study is expected to be of value, in particular, to theorists and engineers that work with nano machines and nano sensing devices at the scale of molecular and thermal fluctuations.
	Specifically, in our work, the Langevin equation
	\begin{subequations}\label{eq:underdamped-Langevin}
		\begin{align}
		\ud X_t &= v_t \ud t\label{eq:underdamped-Langevin-x}\\
		m\ud v_t &= \!-\!\nabla_x U (t,X_t)\ud t \!-\!  \gamma v_t \ud t \!+\! \sqrt{2\gamma k_BT(t)}\ud B_t,\label{eq:underdamped-Langevin-p}
		\end{align}
	\end{subequations}
	provides a model for a molecular system that interacts with a thermal environment. In this,
	$X_t\in \mR$ denotes the location of a particle (which is taken to be one-dimensional for notational simplicity), $v_t\in \mR$ represents its velocity at time $t$, $m$ its mass, $\gamma$ the viscosity coefficient of an ambient medium, $k_B$ the Boltzmann constant, $T(t)$ the time-varying temperature of the heat bath at time $t$ that a statistical ensemble of particles is in contact with, $B_t$ a standard Brownian motion that represents thermal excitation from a heat bath, and $U(t,x)$ a time-varying potential exerting a force $-\nabla_x U (t,x)$ on a particle at location $x\in \mR$. This potential represents a control action that models interaction of our ensemble of particles with the external world. 
	
	This Langevin model has been extensively used to describe heat engines in a Carnot-like cycle, that cycles through contact with heat baths of different temperature, in order to quantify maximum power that can be drawn as well as efficiency at maximum power. 
	Specifically, foundational work has been carried out in the overdamped regime \cite{overdampedseifert,fu2020maximal, VandB2005effic} where inertial effects are neglected,
	as well as in the underdamped regime \cite{underdampedlutz,optimization2017masters} under a low friction assumption, where inertial effects now dominate dissipation. Most importantly, experimental realization and verification of findings have been reported in different settings \cite{Valentin2011exp,Lutz2012exp,Martinez2015exp}. 
	
    Energy and entropy exchange between a heat bath and a particle ensemble is not particular to (thermodynamic) engines. Fluctuations in chemical concentrations in conjunction with the variability of electrochemical potentials may provide a universal source of cellular energy \cite{cilek2009earth}. Indeed, a similar mechanism appears to be at play in biological engines. In these, energy exchange appears to be mediated by continuous processes and energy differentials \cite{yanagida2007brownian}, in contrast to the classical Carnot cycle which alternates between heat baths of constant temperature in an idealized engine. Thus, as contemplated by E.\ Schr\"odinger \cite{schrodinger1944physical,ornes2017core}, the ability of life to maintain a distance from equilibrium may be tied to mechanisms that allow drawing energy out of naturally or self-induced fluctuating temperature or chemical concentrations.

     It is partly the above circle of ideas that led us to consider the rudimentary model of a cyclic process (periodical temperature profile of a heat bath) that fuels a molecular engine represented by stochastically driven Langevin dynamics. Another paradigm that may provide context is that of molecular engines that may some day be engineered to tap onto cyclic temperature fluctuations for work production.
     The work we present below is an extension of 
	recent studies for periodic temperature profile in the linear regime \cite{bauer2016optimal,brandner2015periodic,Bradner2020geom}.
	It addresses the nonlinear regime and, in particular, the case of underdamped dynamics where the low friction assumption applies.

\section{Energetics and low friction approximation}
Throughout the paper we will consider the Langevin dynamics~\eqref{eq:underdamped-Langevin} with a harmonic potential
\[
U(t,x)=\frac{1}{2}q(t)x^2,
\]
where $q(t)$ represents its intensity and constitutes our control variable. 
The total energy of a particle is the sum of kinetic and potential energies, 
\begin{equation*}
E_t = \frac{1}{2}mv_t^2 + \frac{1}{2}q(t)X_t^2.
\end{equation*}
The work and heat exchange  at the level of a single particle, $\dbar W_t$ and $\dbar  Q_t$, respectively, are~\cite{sekimoto2010stochastic}
\begin{subequations}
\label{eq:W-Q}
\begin{align}
\dbar W_t &= \frac{1}{2}\dot{q}(t) X_t^2 \ud t,\\
\dbar  Q_t &= -\gamma v_t^2 {\magenta{\ud t}} + \frac{\gamma k_B T(t)}{m} \ud t + \sqrt{2\gamma k_BT(t)} v_t \ud B_t.
\end{align}
\end{subequations}
These definitions are in agreement with the first law of thermodynamics, i.e. $\ud E = \dbar Q + \dbar W$. 
Taking the expectation of~\eqref{eq:W-Q} yields the average rate of work and heat exchange, 
\begin{subequations}
\label{eq:E-W-Q}
\begin{align}
\dbar \mathcal W_t &= \frac{1}{2}\dot{q}(t) \Sigma_x(t)\ud t,\\ \label{eq:dQ}
\dbar \mathcal Q_t &= \bigg(-\gamma \Sigma_v(t)  + \frac{\gamma k_B T(t)}{m}\bigg)\ud t, 
\end{align}
\end{subequations}
where $\Sigma_x$ and $\Sigma_v$ are the variances of $X_t$ and $v_t$, respectively. The averaged expressions also satisfy the first law of thermodynamics $\ud \mathcal E_t = \dbar \mathcal Q_t + \dbar \mathcal W_t$ where 
\begin{equation*}
\mathcal E_t = \frac{1}{2}q(t)\Sigma_x(t) + \frac{1}{2}m\Sigma_v(t)
\end{equation*} 
is the average total energy of the particle.  

The governing equations for the variance are obtained from the Langevin dynamics~\eqref{eq:underdamped-Langevin}, 
\begin{subequations}\label{eq:lyap_full}
\begin{align}
\dot{\Sigma}_x(t) &= 2\Sigma_{xv}(t),\label{eq:lyap_x}\\
\dot{\Sigma}_{xv}(t) &= \Sigma_v(t) - \frac{q(t)}{m}\Sigma_x(t) - \frac{\gamma}{m}{\Sigma}_{xv}(t),\label{eq:lyap_xv}\\
\dot{\Sigma}_v(t) &= -\frac{2q(t)}{m}{\Sigma}_{xv}(t) - \frac{2\gamma}{m}{\Sigma_v}(t) + \frac{2\gamma k_B}{m^2} T(t),\label{eq:lyap_v}
\end{align}
\end{subequations}
where $\Sigma_{xv}(t)$ is the correlation between $X_t$ and $v_t$. 

Although the differential equations for the evolution of the variance are explicit and in closed-form, the analysis for maximizing efficiency is still challenging due to the presence of the three coupled ODEs. As a result, it is common to consider simplifying approximations that hold in different physical regimes.  

Our analysis concerns the underdamped (or low friction) regime that is studied in   \cite{underdampedlutz,optimization2017masters}. When the temperature $T(t)$ and the harmonic  potential intensity $q(t)$ are constant, the low friction regime holds when the particle in~\eqref{eq:underdamped-Langevin} undergoes many oscillations before dissipation effects become significant. Mathematically, this amounts  to $\sqrt{\frac{q}{m}}\gg \frac{\gamma}{m}$.
{In the case where $T(t)$ and $q(t)$ are time varying periodic functions, it is required in addition that their power spectrum mostly lies below the natural frequency of oscillations $\sqrt{q/m}$. }

Under the low friction assumptions, the equipartition condition \begin{equation}\label{eq:equip}
    \frac{1}{2}q(t)\Sigma_{x}(t)=\frac{1}{2}m\Sigma_{v}(t)
\end{equation} holds approximately \cite[Section 4.2]{optimization2017masters}. 
As a result, it is possible to express the  dynamics of our system as a function of $\Sigma_v(t)$ uniquely.
More precisely, we can use \eqref{eq:lyap_x} to express $\Sigma_{xv}(t)$ as follows:
\begin{equation*}
    \Sigma_{xv}(t)=\frac{1}{2}\dot{\Sigma}_{x}=\frac{m}{2q(t)}\dot{\Sigma}_{v}-\frac{m\dot{q}(t)}{2q(t)^2}{\Sigma_v},
\end{equation*}
where the second equality comes from the equipartition condition ~\eqref{eq:equip}. Now, we can rewrite $\Sigma_{xv}(t)$ in \eqref{eq:lyap_v}, to obtain
\begin{equation}\label{eq:low-fric-Lyapunov}
	\dot{\Sigma}_{v}(t)= \Sigma_{v}(t)\bigg(\frac{\dot{q}(t)}{2q(t)}-\frac{\gamma}{m}\bigg)+\frac{\gamma k_B}{m^2}T(t).\end{equation}
{For this equation to be meaningful, we require $q(t)$ to be (weakly) differentiable.}  Equation~\eqref{eq:low-fric-Lyapunov} is the underlying model for our analysis through the  rest of the paper. 

		\section{Analysis of maximum power}\label{sec:power}
	We consider a stochastic thermodynamic engine driven by a periodic temperature $T(t)$ with period $t_f$. The averaged power output over a cycle is equal to
		\begin{align*}
		   \mathcal P = -\frac{1}{t_f}\int_0^{t_f} \dbar \mathcal W_t. 
		\end{align*}
		Using the first law of thermodynamics  $\ud \mathcal E_t = \dbar \mathcal Q_t + \dbar \mathcal W_t$, the power can be expressed in terms of  heat exchange as follows,
		\begin{align}\nonumber
		   \mathcal P &= \frac{1}{t_f}\int_0^{t_f} (\dbar \mathcal Q_t -\ud  \mathcal E_t)\\
		   &=\frac{\gamma}{t_f} \int_{0}^{t_f}\bigg( \frac{k_B T(t)}{m}-\Sigma_v(t)\bigg)\ud t +  \frac{1}{t_f}(\mathcal E(t_f) -  \mathcal E(0)),\label{eq:power2}
		\end{align}	
		where the definition~\eqref{eq:dQ} is used
		and $\mathcal{E}_t\equiv \mathcal{E}(t)$, for clarity. 
		To account for possible discontinuities in $T(t)$ and $q(t)$, as it is often the case in the underdamped setting~\cite{Lutz2015optic,underdampedlutz,gomez2008under}, we consider a cycle from $0^+$ to $t_f^+$. Thus, $q,\Sigma_v$ and $\mathcal E$, which are periodic, satisfy
		 \begin{equation} \label{eq:periodic-constr} q(0^+)=q(t_f^+) \quad \mbox{and} \quad ~\Sigma_v(0^+) = \Sigma_v(t_f^+), \end{equation}
		 and $\mathcal{E}(0^+) = \mathcal{E}(t_f^+)$.
		 The latter condition removes the boundary terms in \eqref{eq:power2}, while
      conditions \eqref{eq:periodic-constr} impose an integral constraint on $\Sigma_v(t)$. Specifically, if we divide~\eqref{eq:low-fric-Lyapunov} by $\Sigma_v(t)$ and  integrate from $0^+$ to $t_f^+$ we obtain
        \begin{equation}
        \label{eq:integration-lyap}
            \log\left(\frac{\Sigma_v(t_f^+)}{\Sigma_v(0^+)}\right) \hspace*{-3pt}  = \log\left(\frac{q(t_f^+)}{q(0^+)}\right)^{\frac{1}{2}} \hspace*{-5pt} -\frac{\gamma}{m}t_f + \frac{\gamma k_B}{m^2} \int_{0}^{t_f}\hspace*{-3pt} \frac{T(t)}{\Sigma_v(t)}\ud t .
        \end{equation}
      Thus, in view of \eqref{eq:periodic-constr},
     \begin{equation}\label{eq:constraint}
         \int_{0}^{t_f}\frac{T(t)}{\Sigma_v(t)}\ud t  = \frac{ m t_f}{k_B}.   
        \end{equation}		
        This integral constraint is the main ingredient in our analysis for both maximum power and maximum efficiency.
        
		We first consider maximizing power output, namely,
		\begin{align}\label{eq:low-fric-optimization}
&\max_{\Sigma_v(\cdot)} \frac{\gamma}{t_f} \int_{0}^{t_f}\bigg( \frac{k_B T(t)}{m}-\Sigma_v(t)\bigg)\ud t,
	\end{align}
subject to the constraint \eqref{eq:constraint}. We observe that due to \eqref{eq:constraint}, 
        \begin{align*}
            \int_{0}^{t_f} \Sigma_v(t) \ud t  &=   \int_{0}^{t_f} \left(\Sigma_v(t) + c\frac{T(t)}{\Sigma_v(t)}\right) \ud t -   c\frac{ m t_f}{k_B}
        \end{align*}
        for any arbitrary constant $c$.
        The minimal value for the integrand $\Sigma_v(t) +c \frac{T(t)}{\Sigma_v(t)}$, pointwise and with respect to $\Sigma_v(t)$, is attained
        for $\Sigma_v(t) = \sqrt{cT(t)}$. This choice for $\Sigma_v(t)$ also satisfies \eqref{eq:constraint} for $\sqrt c =  \frac{k_B}{mt_f}\int_0^{t_f}\sqrt{T(t)}\ud t$. Therefore,
        		\begin{equation}\label{eq:Sigma-opt}
		 \Sigma_v(t) = \left(\frac{k_B}{mt_f} \int_0^{t_f} \sqrt{T(s)} \ud s\right) \sqrt{T(t)}
		\end{equation} 
       is the unique maximizer of \eqref{eq:low-fric-optimization}.

The relation between $\Sigma_v(t)$ and the protocol $q(t)$ can be obtained from 
 \eqref{eq:low-fric-Lyapunov} by direct integration (cf.\ \eqref{eq:integration-lyap}), which for the optimizing $\Sigma_v(\cdot)$
 gives
		\begin{equation}\label{eq:opt-protocol}
		    q(t) =  q_0 \frac{T(t)}{T(0)} \exp\left( \frac{2\gamma t_f}{m} \left(\frac{t}{t_f} - \frac{\int_0^t\sqrt{T(s)}\ud s}{\int_0^{t_f} \sqrt{T(s)}\ud s}\right)  \right).
		\end{equation}
		Note that, $q_0$ specifies the starting value $q(0)$ of the optimal protocol. 
		Further, $q(t)$ is continuous at all times when $T(t)$ is, which agrees with the intuition that jumps in the control
		are necessitated for adapting to a sudden jump in temperature \cite{overdampedseifert,underdampedlutz}.
		
		Finally,
		inserting \eqref{eq:Sigma-opt} into the expression for power, we obtain that, under this protocol, any value for $q_0$ gives the same level of power, namely,
		\begin{equation}\label{eq:max-P}
		\setlength{\fboxsep}{6pt}
		    \mbox{\fbox{$\calP^* = \frac{\gamma k_B}{m}\text{Var}(\sqrt{T})$}}
		\end{equation}   
		where
			\begin{equation*}
 	  \text{Var}(\sqrt{T})=\frac{1}{t_f}\int_0^{t_f}T(t) \ud t -\left( \frac{1}{t_f}\int_0^{t_f}\sqrt{T(t)}\ud t \right)^2\;
				\end{equation*}
		quantifies fluctuations in $\magenta{\sqrt{T(t)}}$. 

		Thus, expression \eqref{eq:max-P} quantifies the maximum power that can be drawn in terms of an explicit measure of temperature fluctuations. 
		
	It is interesting to consider the shape of temperature fluctuations that allows maximum power to be drawn under the optimal protocol. In view of the above, it is {\purple easy} to see that when $T(t) \in [T_c,T_h]$, the Carnot temperature profile
\begin{align}\label{eq:temp-Carnot}
	T_{\rm Carnot}(t) = \begin{cases}
	T_h,\quad t \in[0,\frac{t_f}{2})\\
	T_c, \quad t \in [\frac{t_f}{2},t_f),
	\end{cases}
	\end{align}		
allows drawing	maximum power, which in this case is
		\[ 
		   \frac{\gamma k_B}{m}\frac{(\sqrt{T_h}-\sqrt{T_c})^2}{4}.
		\] 
	This special profile (Carnot) has already been studied in~\cite{underdampedlutz}. 
	\black
	\section{Analysis of maximum efficiency at fixed power}\label{sec:efficiency}	
The classical definition of efficiency is the ratio of work produced over the amount of heat drawn from the hot bath into the system,
\begin{equation}
\etaQ = \frac{ \purple - \mathcal W}{\;\mathcal Q_h}.
\end{equation} 
This definition presumes that the environment that constitutes the hot heat bath is well defined, and as a consequence that $\mathcal Q_h$ is well defined as well. This is clearly the case
for a Carnot cycle, $T(t)=T_{\rm Carnot}(t)$,
when the system remains in contact with the bath at $T_h$ during the interval $[0,\frac{t_f}{2})$.

For this case, at maximum power,
\begin{equation*}
    \mathcal Q_h = \int_0^{\frac{t_f}{2}} {\red \dbar} \mathcal Q_t = \frac{\gamma k_B}{m}\frac{t_f}{4}\sqrt{T_h}(\sqrt{T_h}-\sqrt{T_c}),
\end{equation*}
where the second equality is obtained using~\eqref{eq:dQ} and the optimal expression for $\Sigma_v$ in~\eqref{eq:Sigma-opt}. 
It follows that
\[ 
    \etaQ = 1- \sqrt{\frac{T_c}{T_h}},
\] 
which is precisely the  Curzon-Ahlborn efficiency \cite{Novikov,CA1975effic}. This result was obtained in \cite{underdampedlutz}.

However, in the setting of an arbitrary periodic temperature profile, it is not entirely clear what constitutes the hot bath. The same bath
serves as hot and cold, at different times, and therefore $\mathcal Q_h$ and $\etaQ$ may be given different definitions and be subject to different interpretations \cite{bauer2016optimal,Bradner2020geom}.

In the present work, following ~\cite{Bradner2020geom}, we adopt
\begin{equation*}
		  \etaU
		  =\frac{\purple -\mathcal{W}}{\;\mathcal{U}},
		\end{equation*}
		as our definition of efficiency, where
		\begin{equation*}
		    \mathcal{U}=-k_B\int^{t_f}_0 S(\rho_t)\dot{T}(t) \ud t,
		\end{equation*}
(also, $\mathcal U=k_B \int T(t) \dot{S}(\rho_t)\ud t$, using integration by parts) represents the effective uptake of thermal energy, while
 \[
 S(\rho_t)=-\int  \int \rho_t(x,v)\log( \rho_t(x,v)) \ud x \ud v
 \]
 denotes the entropy
	of the particle-distribution $\rho_t$. 
		The value $\mathcal U$ can be thought of as the maximal amount of work that can be extracted when the dissipation is zero and the transition takes place quasistatically. In general, over a cycle,
		\[
		\mathcal U= -\mathcal W +\mathcal W_{\rm diss},
		\]
		where the dissipation
		\[
		\mathcal W_{\rm diss} = \int  k_BT(t) \underbrace{\left(\dot{S}(\rho_t) d t - \frac{ \dbar \mathcal Q}{k_BT(t)}\right) }_{\dot S_{\rm total}}, 
		\]
		  relates to the  total entropy production rate $\dot S_{\rm total}$ in the system and the environment.  
For quasi-static transitions $\mathcal W_{\rm diss}=0$ and $\etaU=1$.

	In the underdamped low friction regime with quadratic potential, the distribution $\rho_t$ is Gaussian with independent components $(x,v)$ whose covariance are $\Sigma_x(t)$ and $\Sigma_v(t)$, respectively. Using the standard formula for the entropy of Gaussian distributions together with the equipartition condition ~\eqref{eq:equip} expressed in terms of the variances, $\mathcal{U}$ simplifies to
	\begin{equation*}
		  \mathcal{U}=-\frac{k_B}{2}\int^{t_f}_0 \ln\Big((2\pi e)^2\frac{m}{q(t)}\Sigma_v^2(t)\Big)\dot{T}(t) \ud t.	\end{equation*} 
In case $T(t)$ is discontinuous, the limits of integration need to be taken as $0^+$ and $t_f^+$, and similarly in what follows.
Integration by parts gives
\begin{align*}
    \mathcal{U}=&\frac{k_B}{2}\int^{t_f}_{0}\bigg(2\frac{\dot{\Sigma}_v(t)}{\Sigma_v(t)}-{ \frac{\dot{q}(t)}{q(t)}}\bigg) {\purple T(t) }\ud t\\
    =&\frac{k_B^2\gamma}{m^2}\int^{t_f}_{0}\frac{T(t)^2}{\Sigma_v(t)}\ud t-\frac{k_B\gamma}{m}\int^{t_f}_{0} T(t)\ud t
\end{align*}
where we have used \eqref{eq:low-fric-Lyapunov} to establish the second equality.

We now seek to optimize efficiency at a given power. To this end, it suffices to minimize the dissipation $\mathcal{U}$ subject to the given power. Thus, we consider
\begin{align}\label{eq:low-fric-effic-optimization}
\min_{\Sigma_v(\cdot)}~&\int^{t_f}_{0}\frac{T(t)^2}{\Sigma_v(t)}\ud t,
\end{align}
subject to \eqref{eq:constraint}
together with
\[
\frac{\gamma}{t_f} \int_{0}^{t_f}\bigg( \frac{k_B T(t)}{m}-\Sigma_v(t)\bigg)\ud t = \mathcal{P},
\]
for any given level of power $\mathcal P$.
We follow a variational approach. The Lagrangian for our optimization problem is
\begin{equation*}
    \mathcal{L}=\frac{T(t)^2}{\Sigma_v(t)}+\lambda\bigg(\Sigma_v(t)-\frac{ k_BT(t)}{m}+\frac{\mathcal{P}}{\gamma}\bigg)+\mu\bigg(\frac{T(t)}{\Sigma_v(t)} -\frac{ m }{k_B}\bigg),
\end{equation*}
where $\lambda,\mu$ are the Lagrange multipliers.
The stationarity condition for $\Sigma_v(t)$ gives
\begin{equation*}
   \Sigma_v(t)= \frac{1}{\sqrt{\lambda}}\sqrt{T(t)^2+\mu T(t)}.
\end{equation*}
Applying the power constraint we obtain
\begin{equation*}
    \sqrt{\lambda}=\frac{\int^{t_f}_{0}\sqrt{T(t)^2+\mu T(t)} \ud t}{\frac{ k_B}{m}\int^{t_f}_{0}T(t) \ud t -\frac{t_f}{\gamma}\mathcal{P}}.
\end{equation*}
In a similar manner as before, the optimal protocol can now be obtained by direct integration of \eqref{eq:low-fric-Lyapunov}, to give
	\begin{align}\label{eq:protocol-effic}
		    q(t) =  &q_0 \frac{T(t)^2+\mu T(t)}{T(0)^2+\mu T(0)}  \exp\left( \frac{2\gamma(t - r(t))}{m}   \right),
		\end{align}
where \begin{equation*}
    r(t)=\frac{\int_0^{t_f}\sqrt{T(s)^2+\mu T(s)}\ud s}{\int_0^{t_f} T(s)\ud s-\frac{t_f m}{\gamma k_B}\mathcal{P}}\int^{t}_{0}\frac{\sqrt{T(s)}}{\sqrt{T(s)+\mu}}\ud s.
\end{equation*}
It now remains to fix $\mu$ by imposing the last constraint, which yields $r(t_f)=t_f$. This is,
\begin{align}\label{eq:beta}
    \int^{t_f}_{0}\sqrt{T(t)^2+\mu T(t)}&\ud t \int^{t_f}_{0}\frac{\sqrt{T(t)}}{\sqrt{T(t)+\mu}}\ud t \nonumber\\  =&t_f\int^{t_f}_{0}T(t)\ud t-\frac{t_f^2m}{\gamma k_B}\mathcal{P}.
\end{align}

In general, equation \eqref{eq:beta} cannot be solved explicitly for $\mu$. However, 
it is easily seen that when $\mathcal{P}=0$, then $\mu=0$ is a solution. Also, as 
$\mathcal{P}\to \mathcal{P^*}$, then the solution $\mu\to\infty$. 
It is worthwhile to point out that the condition $r(t_f)=t_f$ makes the protocol \eqref{eq:protocol-effic} continuous for continuous temperature profiles.

\begin{figure}[t!]\includegraphics[width=0.495\textwidth]{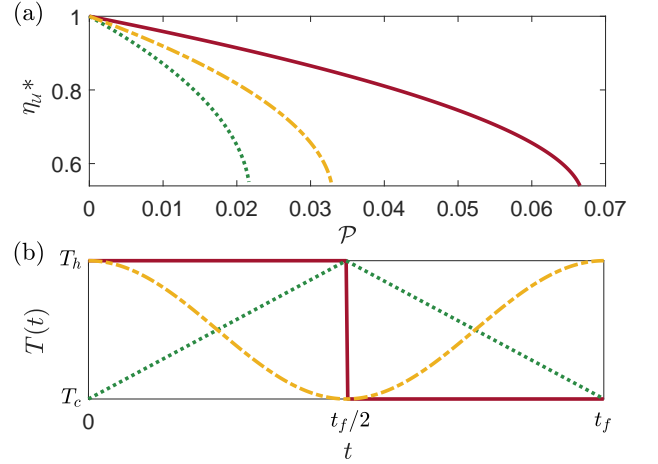}
\caption{ (a) Maximum efficiency at fixed power for the temperature profiles portrayed in (b). (b) Temperature profiles with colors and patterns matching the corresponding performance curves in (a).}
\label{fig:3temp}
\end{figure}

{\red For notational convenience, we let
\begin{equation*}
\overline{g(t)}=\frac{1}{t_f}\int_0^{t_f} g(t)dt
\end{equation*}
denote the average {\purple of a periodic function $g(t)$} over the period.
Now, putting all of the ingredients together, the optimal $\mathcal{U}$ is
\begin{align*}
    \mathcal{U^*}=\frac{t_f}{\kappa}\left[\frac{
    \overline{\sqrt{(T(t)+\mu) T(t)}}}{\overline{T(t)} -\kappa\mathcal{P}}\times \overline{\left(\frac{T(t)^{\frac{3}{2}}}{\sqrt{(T(t)+\mu)}}\right)}-\overline{T(t)}\right]
\end{align*}
where ${\kappa}=m/(\gamma k_B)$. Consequently, the maximum efficiency at fixed power is
\begin{equation}\label{eq:optimal-effic}
  \etaU^*(\mathcal{P})=\frac{\kappa\mathcal{P}}{\left[\frac{
    \overline{\sqrt{(T(t)+\mu) T(t)}}}{\overline{T(t)} -\kappa\mathcal{P}}\times \overline{\left(\frac{T(t)^{\frac{3}{2}}}{\sqrt{(T(t)+\mu)}}\right)}-\overline{T(t)}\right]}.
\end{equation}
}

Figure \ref{fig:3temp} displays the trade-off between power and efficiency for three different temperature profiles.
Note that for vanishing power we obtain the maximum possible efficiency, i.e. $\etaU^*=1$, that corresponds to vanishing dissipation. On the other hand, as {\red $\kappa \mathcal{P}\to \kappa \mathcal{P}^* = \text{Var}(\sqrt{T})$ (from \eqref{eq:max-P})}, the limiting value from \eqref{eq:optimal-effic} gives that
\begin{equation}\label{eq:effic_maxp}
   \etaU^*\to \frac{
  \overline{\sqrt{T(t)}}
  \text{Var}(\sqrt{T}) }
    {\overline{ T(t)^\frac{3}{2}} -
    \overline{T(t)}\times \overline{\sqrt{T(t)}}
    }
    .
\end{equation}
This is precisely the efficiency at maximum power  {\red obtained by using the optimal protocol \eqref{eq:opt-protocol}}. {The limit in \eqref{eq:effic_maxp} can also be expressed in terms of the third moment of the square root of the temperature 
\begin{equation*}
    \mu_3(\sqrt{T})=\frac{1}{\; t_f}\int_0^{t_f}\Big(\sqrt{T(t)}-\overline{\sqrt{T(t)}}\Big)^3\ud t, \end{equation*}
    specifically,
\begin{equation}
 \setlength{\fboxsep}{8pt}
 \mbox{\fbox{
  $\etaU^* \to$ \begingroup \Large$\frac{1}{  2+ \frac{   \mu_3(\sqrt{T})   }{    \text{Var}(\sqrt{T})\overline{\sqrt{T(t)}}   } }.$\endgroup
  }}
\end{equation}
Thus, for any  temperature profile such that $\mu_3(\sqrt{T})=0$, the limit is $\frac12$, which in particular is true for the Carnot-like temperature profile $T(t)=T_{\rm Carnot}(t)$. We also note that the third moment can be negative, leading to efficiency that exceeds $\frac12$, albeit at the cost of decreasing maximum power.}

\section{Example: Sinusoidal temperature profile}
In the following we compare our results to the ones obtained in the linear response regime, where it is assumed that $T_h-T_c\ll T_h+T_c$ \cite{bauer2016optimal,Fu2021Harvesting}. To this end, we consider 
		the sinusoidal temperature profile
    \begin{equation}
    T(t)=\frac{T_h+T_c}{2}+\frac{T_h-T_c}{2} \cos(\omega t),\label{eq:tempprofile}
    \end{equation}
		with $\omega\ll\sqrt{\frac{ q}{m}}$ so as to ensure that the low friction assumption holds. 
		
			We first use our analytical result to recover the expressions obtained in the linear response regime. Specifically, assuming $\Delta T=\frac{T_h-T_c}{2}$ is much smaller than  $\bar{T}=\frac{T_h+T_c}{2}$, we have that
			\begin{align*}
		    q^*(t)&=q_0+q_0\frac{ \Delta T}{\bar{T}}\cos(\omega t) +\mathcal{O}(\Delta T^2),\\
		    \cal{P}^*&=\frac{\gamma k_B }{8 m }\frac{\Delta T^2}{\bar{T}}+\mathcal{O}(\Delta T^3),\\
		    \etaU^*(\mathcal{P^*})&=\frac{1}{2}+\mathcal{O}(\Delta T^2),
	\end{align*}
		where the first two coincide with the expressions for the optimal protocol and power that are given in~\cite{bauer2016optimal,Fu2021Harvesting}, in the limit where  $\frac{q_0}{m}$ is large enough to satisfy the requirements of the low friction regime. In fact, in the linear response regime
$\etaU^*(\mathcal{P^*})=\frac{1}{2}+\mathcal{O}(\Delta T^2)$ holds for any temperature profile whose time dependent component is an odd function of $t$. 
		
		 \begin{figure}[t!]
\includegraphics[width=0.45\textwidth]{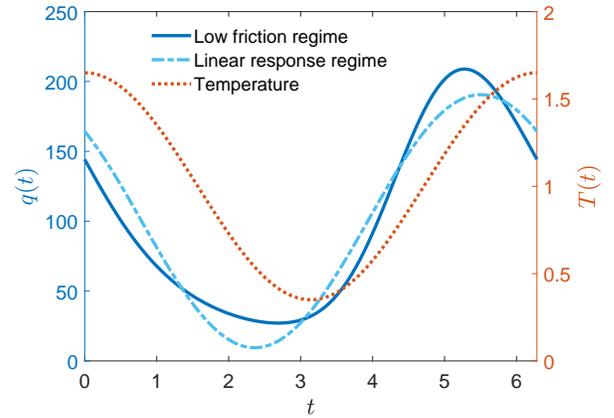}
\caption{Comparison between the optimal protocols $q(t)$ in the low friction and linear response regimes, for the sinusoidal temperature profile shown with a red-dotted curve.}
\label{fig:lowlin_qT}
\end{figure}
\begin{figure}[t] \centering
\includegraphics[width=0.51\textwidth,trim={0pt 0pt 0pt 0pt}]{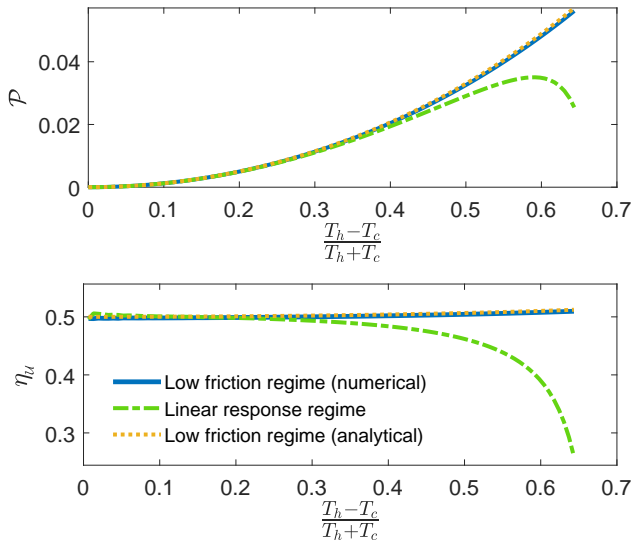}
\caption{{Power and efficiency for \eqref{eq:underdamped-Langevin} with a temperature profile \eqref{eq:tempprofile} as functions of temperature  ratios for two protocols:
i) the optimal protocol obtained under the low friction assumption~\eqref{eq:opt-protocol} (solid line), and
ii) the optimal protocol obtained in the linear response regime  \cite[Eqs 29, 30a]{Fu2021Harvesting} (dashed line).}
For comparison, the analytical expressions for maximum power~\eqref{eq:max-P} and efficiency~\eqref{eq:effic_maxp} in the low friction regime are also drawn (dotted curve).}
 \label{fig:lowlin_pe}
\end{figure}

{Moreover, in order to illustrate the utility of our results, we compare the optimality of the protocol obtained in this paper, to the one obtained in the linear response regime. Specifically,} we solve numerically the equations for the variance \eqref{eq:lyap_full} for two protocols:
\begin{itemize}
\item[i)] the optimal protocol obtained with low friction assumption according to~\eqref{eq:opt-protocol},
\item[ii)] the optimal protocol obtained in linear response regime in \cite{bauer2016optimal} and~\cite[Eqs 29, 30a]{Fu2021Harvesting}.
\end{itemize}
These two protocols, along with the temperature profile, are shown in Figure~\ref{fig:lowlin_qT}. Then, we use the numerical solution to the equations for the variance to compute power and efficiency. The result, as the temperature difference $T_h-T_c$ is varied, is shown in Figure~\ref{fig:lowlin_pe}. 

In these numerical examples, we made sure  $\frac{\gamma}{m}\ll \sqrt{\frac{q}{m}}$  so that the low friction assumption remains valid. { As it is expected, the optimal protocol under the linear response regime assumption performs well when the temperature ratio is small, giving a numerical result that agrees with equations (13) and (20). However, as the temperature ratio increases, higher order terms become significant and the protocol is no longer optimal.}

\begin{figure}[h]\centering 
\includegraphics[width=0.5\textwidth]{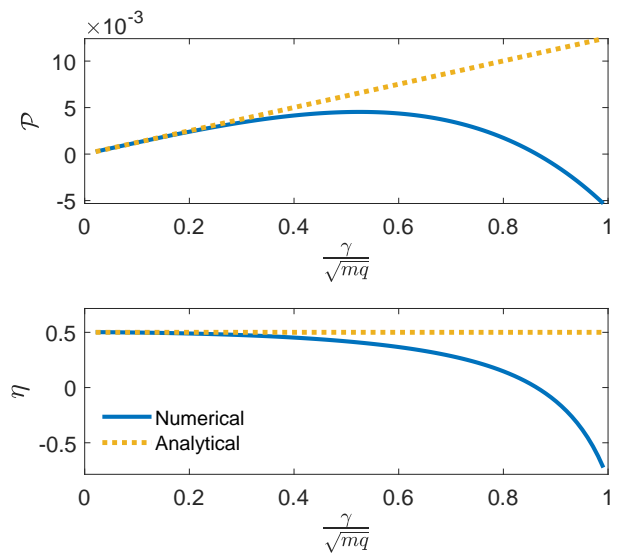}
\caption{
Numerical values for power and efficiency as a function of the friction coefficient for the system governed by equation \eqref{eq:lyap_full}, driven by the optimal protocol~\eqref{eq:opt-protocol} obtained in the low friction regime (solid line). For  comparison, we have included the analytical expressions for maximum power~\eqref{eq:max-P} and the efficiency at maximum power~\eqref{eq:effic_maxp}  obtained under the low friction assumption (dotted line).}
\label{fig:degradation_with_gamma}
\end{figure}

{	
We finally highlight how the numerical value of power and efficiency, using the optimal protocol \eqref{eq:opt-protocol}, deviates from the analytical 
expressions for maximum power~\eqref{eq:max-P} and efficiency at maximum power~\eqref{eq:effic_maxp} with increasing values of the friction coefficient $\gamma$. This is shown in Figure \ref{fig:degradation_with_gamma}. Specifically, for small values of $\frac{\gamma}{\sqrt{mq}}$, the analytical and numerical results for power and efficiency coincide. For larger values of $\frac{\gamma}{\sqrt{mq}}$ the approximation error becomes significant due to the fact that the condition $\frac{\gamma}{\sqrt{mq}}\ll 1$, and hence, the approximate equipartition relation ~\eqref{eq:equip} , no longer hold. 
}

\section{Concluding remarks}\label{sec:conclusion}
    	{\red 
    	The present work quantifies maximum power, as well as the efficiency at any level of power, 
    	that can be achieved by a thermodynamic engine in the (underdamped) low friction regime, modeled by a Langevin equation and in contact with a heat bath with periodic temperature profile.
    	
    	It turns out that the maximum power and efficiency are expressed in terms of finite moments of the square root of the temperature profile. Moreover, when the temperature changes continuously, the optimal control protocol is expressed explicitly in terms of the temperature and certain finite moments, and is continuous in time as well.
    	
    	Of particular interest is the efficiency at maximum power and, specifically, the relation between useful work and dissipation. 
    	Interestingly, when the third moment of the square root of the temperature profile is zero (for instance, due to a suitable time symmetry), at maximum power in the low {\purple friction} regime, the losses in dissipation equal the amount of work that can be extracted by the engine, yielding $\etaU=\frac12$.
    	{\purple It is of interest to explore possible connections of this result to the universal bound on the efficiency at maximum power being $\frac{1}{2}\eta_C$, where $\eta_C$ is the Carnot efficiency~\cite{esposito2009universality,VandB2005effic}.}
    	}

	\section{Acknowledgements} OMM, RF and AT contributed equally and AT oversaw the technical development. The research was partially supported by NSF under grants 1807664, 1839441, 1901599, 1942523, and
AFOSR under FA9550-20-1-0029.

\bibliography{main}

\end{document}